\documentclass[a4paper,10pt]{amsart}
\usepackage[latin1]{inputenc}

\usepackage{amsfonts}
\usepackage{amsmath}
\usepackage{amssymb}
\usepackage{amsthm}

\usepackage{enumerate}
\usepackage{tikz}
\usetikzlibrary{matrix,arrows}

\newcommand{\an}{\text{an}}

\newcommand{\GG}{\mathbb{G}}

\newcommand{\PP}{\mathbb{P}}

\newcommand{\RR}{\mathbb{R}}

\newcommand{\mF}{\mathcal{F}}

\newcommand{\fX}{\mathfrak{X}}

\newcommand{\trop}{\operatorname{trop}}
\newcommand{\im}{\operatorname{im}}
\renewcommand{\div}{\operatorname{div}}
\newcommand{\supp}{\operatorname{supp}}
\newcommand{\PGL}{\operatorname{PGL}}
\newcommand{\skeleton}{\mathcal{S}}
\newcommand{\longhookrightarrow}{\lhook\joinrel\relbar\joinrel\rightarrow}

\newcommand{\<}{\langle}
\renewcommand{\>}{\rangle}
\newcommand{\inv}{^{-1}}

\newcommand{\quotient}[2]{\left.\raisebox{.2em}{$#1$}\middle/\raisebox{-.2em}{$#2$}\right.}

\newcommand{\rip}{r_i}
\newcommand{\rim}{R_i}
\newcommand{\zip}{b_i}%{z_i^+}
\newcommand{\zim}{c_i}%{z_i^-}
\newcommand{\zjp}{b_j}%{z_j^+}
\newcommand{\zjm}{c_j}%{z_j^-}
\newcommand{\zeinsp}{b_1}%{z_1^+}
\newcommand{\zeinsm}{c_1}%{z_1^-}
\newcommand{\reinsp}{r_1}
\newcommand{\reinsm}{R_1} %Radius der Kreisscheibe zu \gamma_1^-1
\newcommand{\zweip}{b_2}%{z_2^+}

\newcommand{\zweim}{c_2}%{z_2^-}

\newcommand{\discB}{B}
\newcommand{\discC}{C}
\newcommand{\cycle}{Z}
\newcommand{\curve}{X}

\newcommand{\gaussCi}{\zeta_i^-} %{\zeta_{\zim, \rim}}
\newcommand{\gaussBi}{\zeta_i^+} %{\zeta_{\zip, \rip}}

\newcommand{\Trop}{\operatorname{Trop}}

\theoremstyle{definition}

\newtheorem*{theorem*}{Theorem}
\newtheorem{theorem}{Theorem}[section] %chapter]
\newcommand{\Thm}[2]{
\begin{theorem}
\label{thm:#1}
#2
\end{theorem}
}

\newtheorem{lemma}[theorem]{Lemma}
\newcommand{\Lemma}[2]{
\begin{lemma}
\label{lemma:#1}
#2
\end{lemma}
}

\newcommand{\Prf}[1]{
\begin{proof}
#1
\end{proof}
}

\newtheorem{proposition}[theorem]{Proposition}
\newcommand{\Prop}[2]{
\begin{proposition}
\label{prop:#1}
#2
\end{proposition}
}

\newtheorem{corollary}[theorem]{Corollary}

\newtheorem{definition}[theorem]{Definition}

\newtheorem{example}[theorem]{Example}

\newtheorem{remark}[theorem]{Remark}

\newcommand{\fig}[3]{
\begin{figure}
\centering
#3
\caption{#2}
\label{fig:#1}
\end{figure}
}

\title{Faithful Tropicalization of Mumford Curves of Genus Two}
\author{Till Wagner}
\address{Goethe Universit\"at Frankfurt \\ Institut f\"ur Mathematik \\ Robert-Mayer-Str. 6-10 \\ 60325 Frankfurt am Main \\ Phone: 069 798 28281 \\ Fax: 069 798 22302}
\email{wagner@math.uni-frankfurt.de}
\date{}
\keywords{Mumford Curves, Tropical Geometry, Berkovich spaces, Tropicalization}
\subjclass[2000]{14T05, 14H30, 14H45}

\begin{document}

\maketitle

\begin{abstract}

In the present paper, we investigate the question if the skeleton of a Mumford curve of genus two can be tropicalized faithfully in dimension three, i.e. if there exists an embedding of the curve in projective three space such that the tropicalization maps the skeleton of the curve isometrically to its image. 
Baker, Payne and Rabinoff showed that the skeleton of every analytic curve can be tropicalized faithfully. However the dimension of the ambient space in their proof can be quite large. 

We will define a map from the skeleton to the tropicalization of the Jacobian, which is an isometry on the cycles. 
It allows us to find principal divisors and simultanously to determine the retractions of their support on the skeleton which is necessary to calculate their tropicalization. 

It turns out that a Mumford curve of genus two whose cycles of its skeleton are either disjoint or share an edge of length at most half of the length of the cycles, can be tropicalized faithfully in dimension three. 
\end{abstract}
\section{Introduction}

This paper is inspired by 
recent work of Baker, Payne and Rabinoff (\cite{bpr} and \cite{bpr2}) which compares the  analytic and tropical geometry of curves.
Let $K$ be a complete, non-archimedean field with a non-trivial absolute value, and 
 let $X$ be a smooth, proper and connected curve over $K$.  Then the subset $H_\circ(X^\an)$ of non-leaves in the Berkovich space $X^{\an}$ carries a natural metric which is induced by skeletons of semistable models. 
 If $X$ is embedded in a toric variety and meets the dense torus $T$, the resulting tropicalization $\Trop(T \cap X)$ is a  one-dimensional polyhedral complex. All its edges  have rational slopes with respect to the cocharacter lattice of the torus. Hence $\Trop(T \cap X)$ can be endowed with a natural metric  locally given by the lattice length on each edge. In general the tropicalization map $(X \cap T)^{\an} \rightarrow \Trop(X \cap T)$ is not injective and does not respect the metrics on both sides. 
 
Among the main results in \cite{bpr2} are Theorem 6.20 and 6.22, which say that for every skeleton $\Gamma$ in  $H_\circ(X^\an)$  there exists a closed immersion of $X$ in a quasiprojective toric variety such that $X$ meets the dense torus $T$ and $\Gamma$ maps isometrically to its image under the induced tropicalization map on $X \cap T$. 
In the terminology of \cite[6.15.2]{bpr2} we say that in this case $\Gamma$ is faithfully tropicalized. 

Given a skeleton of an analytic curve, it is an interesting question to determine the minimal dimension necessary to tropicalize it faithfully. 
Using the construction in the proof of Theorem 6.22 in \cite{bpr2}, we may assume that the skeleton has at least two vertices and there are at least $g+1$ edges. To tropicalize the skeleton faithfully they define two functions for each edge as well as one function for each pair of vertices. These functions may not define an embedding thus we have to add more rational functions (for example four additional functions since every curve can be embedded in projective space of dimension three). This sums up to at least $2g+6$ functions. Given a finite subset $D\subset X(K)$ for every point $p$ in $D$ one needs two functions to tropicalize the ray leading to $p$ faithfully. Hence we obtain an embedding of dimension at least $2(g+|D|)+5$ that leads to a faithful tropicalization of a skeleton. 

However, in special cases one can do much better. For example, for Tate curves, this bound would be 7, whereas 
by \cite[Theorem 7.2]{bpr2} there exists a faithful tropicalization in twodimensional projective space. Furthermore Chan and Sturmfels showed that after a change of coordinates every tropicalization of an elliptic curve is in so called honeycomb form which is a faithful tropicalization (\cite{chst}, Theorem 7).

In \cite[Theorem 8.2]{bara} Baker and Rabinoff defined for each curve $X$ together with a skeleton $\Gamma$ a rational morphism $X\dashrightarrow \PP^3$ whose tropicalization, when restricted to $\Gamma$, is an isometry onto its image. However this morphism may not be an embedding. 

In the present paper, we investigate this question for Mumford curves of genus two. 
There are two types of such curves: the cycles in the skeleton can share an edge or the cycles are disjoint. 
It turns out that for curves whose cycles are disjoint there exists an embedding in projective three-space such that the tropicalization is faithful (Theorem \ref{thm:faithful_extended_skeleton}). This is also true for curves whose cycles share an edge if the cycles are at least twice as long as the shared edge. 

Let $D$ be a very ample divisor of degree five such that its support retracts to given points on the skeleton. We define two divisors as follows: Choose three points from the support of $D$ (three poles) and one additional $k$-rational point on the curve (one zero). From general theory we know that there exist two $k$-rational points (zeros) so that the resulting degree zero divisor is principal, i.e. we obtain a rational function on the curve with three poles. However we do not know where the retractions of these new zeros lie on the skeleton. This is necessary to calculate the slopes of the tropicalization.% of these functions with the slope formula.

To solve this problem we will make use of the tropicalization of the Jacobian, which inherits a group structure from the Jacobian, and a map $\mu$ from the skeleton to the tropicalization of the Jacobian, which is an isometry on the cycles. 

They will tells us where the new zeros of the function are retracting to (Lemma \ref{lemma:points}) and we are able to show that there exists a rational morphism from the curve to twodimensional projective space tropicalizing the skeleton faithfully. 

Next we check if the rays from the extended skeleton are also tropicalized faithfully. We will see that two rays may intersect. To correct this we have to define a third function in the linear system $L(D)$, linearly independent from the others. These three functions tropicalize the extended skeleton faithfully, and, since they are generators of $L(D)$, they also define an embedding of the curve in projective three-space.

\section{Skeletons and tropicalizations of curves}

Let $K$ be a algebraically closed field which is complete with respect to a non-Archimedean non-trivial absolute value $|\cdot|$. As usual, we put $K^\circ = \{ x \in K: |x| \leq 1\}$ and $K^{\circ\circ} = \{ x \in K: |x| < 1\}$ and we denote by $\widetilde{K} = K^\circ / K^{\circ\circ}$ the residue field of $K$.% We fix an algebraically closed and complete extension field $\mathbb{C}_K$ of $K$. 

For every $K$-scheme $X$ of finite type we denote by $X^\an$ the associated Berkovich analytic space, as defined in \cite{ber}. If $X$ is a smooth, projective and geometrically connected curve of positive genus, then $X^\an$ is a special quasipolyhedron in the sense of \cite{ber}, Definition 4.1.1 and 4.1.5. Its Betti number is at most $g$ by \cite{ber}, Theorem 4.3.2.

Let $X$ be a smooth, projective curve over $K$ with semistable reduction, and fix a semistable model  $\fX$ of $X$. This gives rise to a skeleton $\skeleton(\fX)$ inside $X^{\an}$ which is a deformation retract, see \cite[Chapter 4]{ber}. Let $\tau: X^\an\rightarrow \skeleton(\fX)$ be the continuous retraction map .

Let $D$ be a finite subset of $X(K)$ and $V$ be a finite set of type-2 points in  $X^\an$. Set $U=X \setminus D$. Then $U^\an\setminus V$ is a \emph{semistable decomposition} in the sense of \cite[Definition 3.1]{bpr} if it is the disjoint union of
\begin{itemize}
 \item infinitely many open balls,
 \item finitely many open annuli,
 \item finitely many punctured balls, each of them containing precisely one point of $D$.
\end{itemize}

The punctured balls can be seen as open annuli with infinite modulus.
Their skeleton is a semi-infinite ray (see \cite{bpr}, chapter 2.1).

For $D$ as before set $W=\{\tau(x):x\in D\}\cup\{$vertices of $\skeleton(\fX)\}\subset \curve^\an$. 
This is a semistable vertex set of $X^{\an}$  in the sense of \cite{bpr}, Definition 3.1. 
If two distinct points $x,y$ of $D$ lie in the same open ball $B$ of $X^\an\setminus W$ then there exists a smallest closed disc in $B$ containing $x$ and $y$. The Gauss point of this closed disc will be denoted by $x\vee y$. If we adjoin all such points $x\vee y$ to $W$ then it is still a semistable vertex set and $U^\an \setminus W$ is a semistable decomposition. This follows from the proof of \cite[Lemma 3.13 (3)]{bpr}. 
The path from $x\vee y$ to $\skeleton(\fX)$ yields a new bounded edge. 
The skeleton associated to these data in \cite[Definition 3.3]{bpr} is
the union of $\skeleton(\fX)$ with the new bounded edges, the vertices $W$ and the rays leading to points in $D$.
We denote it by $\skeleton(\fX, D)$ and call it an extended skeleton.
By \cite[Lemma 3.8]{bpr} there exists a continuous retraction map $\tau: U^{\an} \rightarrow \skeleton(\fX, D)$ from the open curve to the extended skeleton. 

A skeleton carries a natural metric, the shortest-path metric. The distance between two points is the smallest length of all paths connecting the points, where the length of an edge is defined to be the logarithmic modulus of the corresponding annulus (\cite[Def. 3.10 and chapter 3.12]{bpr}). This metric extends to rays giving the skeleton a degenerate metric. 

We will use the following slope formula (\cite[Theorem 5.15]{bpr}),
see also Thuillier's thesis \cite[3.3.15]{th}.

\begin{theorem}[Slope formula]
\label{thm:poincare}
Let $f$ be a rational function on $X$ and $D$ be a finite subset of $X(K)$ containing $\supp(\div(f))$. Set $U=X\setminus D$, and consider the retraction map $\tau: U^{\an} \rightarrow \skeleton(\fX, D)$. Then the function $\log|f|:U^\an\to\RR$, mapping a point $x$ to $\log|f(x)|$, has the following properties:

\begin{enumerate}[i)]
\item $\log|f|=\log|f|\circ\tau$, i.e. $\log|f|$ factors over the skeleton.
\item $\log|f|$ is piecewise affine-linear with integer slopes on each edge of the skeleton $\skeleton(\fX,D)$. 
\item Let $x$ be a point in the skeleton and let $d_v(x)$ denote the outgoing slope of $\log|f|$ along an edge $v$. Then $\sum_v d_v(x)=0$ where only finitely many outgoing slopes are unequal to zero.
\item If $v$ is a vertex in $\skeleton$ and $r$ a ray in $v$ leading to a point  $x \in \supp(\div(f))$, let $d_v$ be the outgoing slope of $\log|f|$ along $r$. If $x$ has  multiplicity $m$ in $\div(f)$, then $d_v(x)=-m$.
\end{enumerate} 
\end{theorem}

If $x$ is a point in the skeleton that is neither a vertex nor a retraction of a pole or zero then there are two edges of the skeleton emanating from $x$ thus by iii) the slopes along theses edges have to be equal. Hence the slope of $\log|f|$ can only change at the retractions of the poles and zeros of $f$ and the vertices of the skeleton.

Let $U$ be a connected $K$-scheme together with a closed immersion $\varphi: U \hookrightarrow \GG_m^n$ into a $K$-split torus. As a set, the associated tropical variety $\Trop(U)=\Trop_\varphi(U)$ is the image of $U^\an$ under the tropicalization map 
\[(\GG_m^{n})^{\an} = (\mathrm{Spec} K[x_1^\pm, \ldots, x_n^\pm])^\an \longrightarrow \RR^n, \quad p \longmapsto (\log|x_1(p)|, \ldots, \log|x_n(p)|).\]
Note that $\Trop_\varphi(U) = \Trop_{\varphi_L}(U_L)$ for any base change by a non-Archimedean complete extension field $L/K$, see \cite{gu}, Proposition 3.7. So if we start with a non-algebraically closed field we can perform a base change to apply the theory developed in this paper. 
The tropicalization $\Trop(U)$ can be enriched with the structure of a
balanced, weighted, integral polyhedral complex of pure dimension 
$d = {\rm dim}(U)$, see e.g. \cite[Theorem 3.3.6]{mast}. %
$\Trop(U)$ is called faithful if the restriction of the tropicalization map to the skeleton of the curve is an isometry. If we consider extended skeleta we call the tropicalization faithful if it is faithful on the skeleton, homeomorphic on the extended skeleton and the rays have slopes $\pm1$. 
Let $\varphi:X^\an\to(\PP^n)^\an$ be an embedding in projective space given by rational functions $f_0,\dots,f_n$ and let $D=\bigcup_i \supp(\div(f_i))$ be the union of the support of the corresponding principal divisors. We say that the curve is tropicalized faithfully with respect to $\varphi$ if $\Trop(X\setminus D)$ is a faithful tropicalization of the extended skeleton $\skeleton(\fX,D)$.

\section{Mumford curves}

The present paper deals with analytic Mumford curves, for which the Betti number of $X^\an$ is equal to the genus. Mumford curves are curves with totally degenerate stable reductions. In \cite{mu} Mumford shows that these are precisely the curves admitting a non-archimedean Schottky uniformization. Let us recall some facts about these uniformizations.

Let $\Gamma$ be a subgroup of $\PGL(2,K)$. A point $ x \in \mathbb{P}^1(K)$ is called a limit point of $\Gamma$, if there exists some $y \in  \mathbb{P}^1(K)$ and an infinite subset $\{\gamma_n: n \geq 1\}$ of $\Gamma$ such that $\gamma_n(y) \rightarrow x$. We denote by $\mathcal{L} = \mathcal{L}_\Gamma$ the set of all limit points. A subgroup $\Gamma$ of $\PGL(2,K)$ is called a Schottky group if it is finitely generated, free and discontinuous, where the last property means that $\mathcal{L} \neq \mathbb{P}^1(K)$ and that  all orbit closures in $\mathbb{P}^1(K)$ are compact. Every element $\gamma \neq 1$ of a Schottky group is hyperbolic, i.e. it has two different fixed points with different absolute values.

Let $\Gamma$ be a Schottky group, and write $\Omega = (\mathbb{P}^1_K)^\an \backslash \mathcal{L}$, where $(\mathbb{P}^1_K)^\an$ is the Berkovich analytic space associated to the projective line over $K$.  Then $\Gamma$ acts freely and properly discontinuously on $\Omega$, and the quotient $\Omega / \Gamma$ is a proper analytic curve over $K$, and hence the analytification of an algebraic curve $\curve$. Every curve over $K$ with  such a Schottky uniformization is called a Mumford curve.

It is shown in \cite{gevp}, section 1, paragraph 4, that every Schottky group has a good fundamental domain, i.e. for each Schottky group $\Gamma$ of rank $g$ %
%where $\infty$ is not a limit point,
there exist 
$2g$ pairwise disjoint ``closed'' discs $\discB_1, \ldots, \discB_g$, $\discC_1, \ldots, \discC_g$ with centers in $K$ and generators $\gamma_1, \ldots, \gamma_g$ of $\Gamma$ such that 
\begin{itemize}
\item $ \mathcal{F} = \mathbb{P}^1(K) \backslash (\bigcup_{i} \discB_i \cup \bigcup_i \discC_i^\circ)$ is a fundamental domain for the action of $\Gamma$ on $\Omega$
\item $\gamma_i \left(  \mathbb{P}^1(K) \backslash \discB_i^\circ\right) = \discC_i \mbox{ and } \gamma_i\inv\left( \mathbb{P}^1(K) \backslash \discB_i\right) = \discC_i^\circ$ for all $i=1,\dots,g$.
\end{itemize}
We denote the corresponding ``open'' discs by $B_i^\circ,\ C_i^\circ$ etc. 

Note that we always choose our coordinates on $\mathbb{P}^1$ so that $0$ lies outside these discs. %

We fix centers $\zip$ of $\discB_i$ and $\zim$ of $\discC_i$. Let $\rip$ be the radius of $\discB_i$ and $\rim$ the one of $\discC_i$. 
Since $0$ is not contained in any disc, the radii of $\discB_i$ and $\discC_i$ satisfy $\rip < |\zip|$ and $\rim < |\zim|$. Furthermore, the absolute value of all points in $\discB_i$ is equal to $|\zip|$ and the absolute value of all points in $\discC_i$ is equal to $|\zim| $. Since $\discB_i$ contains one fixed point of $\gamma_i$ and $\discC_i$ contains the other, and since the two fixed points of a hyperbolic element have different absolute values, we have $|\zip| \neq |\zim|$. Moreover, for every $a \in \discB_i$ and $b \notin \discB_i$ we have $|b-a| = |b -\zip|$. Similarly, for every $a \in \discC_i$ and $b \notin \discC_i$, we have $|b-a| = |b-\zim|$.

Now choose some $a \in \mF$ and define for all $i = 1, \ldots g$
\[u_i(z)=\prod\limits_{\gamma\in\Gamma}\frac{z-\gamma a}{z-\gamma\gamma_i a}.\]
By \cite{gevp}, section II.2, $u_i(z)$  is an analytic function without zeros on $\Omega$ which is independent of the choice of $a$. Moreover, the function $u_i(z)$ is an automorphic function with respect to $\Gamma$, i.e. the quotient $c_i(\gamma) = u_i(z) / u_i(\gamma z)$ is a constant in $K$ independent of $z$ by \cite{gevp}, chapter VI, paragraph 2. %\cite{madr} Proposition 2.
This constant is multiplicative in $\gamma$ and satisfies $c_i(\gamma_j) = c_j(\gamma_i)$. We put $q_{ij} = c_i(\gamma_j)$. Then the matrix $(q_{ij})_{ij}$ is symmetric.

\Prop{ufunction}{
For every $z \notin \bigcup\limits_{i = 1}^g(\discB_i^\circ \cup \discC_i^\circ)$ we have $|u_i(z)|=\left|\frac{z-\zip}{z-\zim}\right|$. 
}

\Prf{
Choose $a \in \mathbb{P}^1(K)$ inside the complement of the union of the discs $\discB_i, \discC_i$. Let $\gamma$ be in $\Gamma$ such that $\gamma \neq 1$ and $\gamma \neq \gamma_i$. Then $\gamma = \gamma_{i_m}^{s_m} \gamma_{i_{m-1} }^{s_{m-1} } \ldots \gamma_{i_1}^{s_1}$ for $i_1, \ldots, i_m \in \{1, \ldots, g\}$ and signs $s_1, \ldots, s_m \in \{-1, 1\}$ . We assume that this expression is reduced, i.e. that  $\gamma_j $ and $\gamma_j^{-1}$ are never adjacent factors.  Let us show by induction on $m$ that the points $\gamma a$ and $\gamma \gamma_i a$ are contained in the same disc. Since $\gamma_j$ maps the complement of $\discB_j $ to $\discC_j^\circ$  and $\gamma_j^{-1}$ maps the complement of $\discC_j$ to $\discB_j^\circ$, the point $\gamma_ia$ is contained in $\discC_i^\circ$ which lies in the complement of all $\discB_j$ and in the complement of those $\discC_j$ with $j \neq i$. Hence our claim holds for $m = 1$. 

If $\gamma$ is an element of length $m > 1$ in $\Gamma$, consider $\gamma' = \gamma_{i_{m-1} }^{s_{m-1}}  \ldots \gamma_{i_1}^{s_1}$. Then $\gamma' \neq 1$. If $\gamma' = \gamma_i^{-1}$, then $\gamma = \gamma_j^{s} \gamma_i^{-1}$, where $s \in \{-1, 1\}$ such that $\gamma_j^{s } \neq \gamma_i$. In this case, the points $\gamma_i\inv a$ and $a$ lie in the complement of $\discB_j$ (resp. $\discC_j$),  hence they are mapped to the disc $\discC_j^\circ$ (resp. $\discB_j^\circ$) by $\gamma_j^s$, and our claim holds. Therefore we may assume that $\gamma'$ is neither $1$ nor $\gamma_i$, in which case we conclude by induction hypothesis that $\gamma' a$ and 
$\gamma' \gamma_ia $ are contained in $\discC_{i_{m-1}}^\circ$ (resp. $\discB_{i_{m-1}}^\circ$). Since $\gamma_{i_m}^{s_m} \neq \gamma_{i_{m-1}}^{-s_{m-1}}$, the element $\gamma_{i_m}^{s_m}$ maps $\discC_{i_{m-1}}$ to $\discB_{i_m}^\circ$ (resp. $\discB_{i_{m-1}}$ to $\discC_{i_m}^\circ$) and our claim follows. 

Hence for all $\gamma$ not equal to $1$ or $\gamma_i\inv$ we have
$|z-\gamma\gamma_i a| = | z - \gamma a|$, since $z$ lies  outside every open disc. This implies that 
\[ |u_i(z)|=\left|\prod\limits_{\gamma\in\Gamma}\frac{z-\gamma a}{z-\gamma\gamma_i a}\right| 
         =\frac{|z-\gamma_i^{-1}a|}{|z-\gamma_i a| } = \frac{|z- \zip|} {|z-\zim|},\]
         since $\gamma_i^{-1}a$ lies in $\discB_i$ and $\gamma_i a$ lies in $\discC_i$. 
         
         }

Let us recall some facts about the Berkovich projective line $(\mathbb{P}^1)^\an = (\mathbb{A}^1)^\an \cup \{\infty\}$. 
The points in the affine line $(\mathbb{A}^1)^\an$ correspond to muliplicative seminorms on the polynomial ring $K[T]$ extending the absolute value on $K$. 
These seminorms can be described explicitly and are classified into four types, see \cite{ber}, 1.4.4.
Every point $a$ in $\mathbb{A}^1(K)$ induces a point of type one, i.e. the multiplicative seminorm
$ f \mapsto |f(a)|$.  For every point $a \in \mathbb{A}^1(K)$ and every positive real number $r$, there is a corresponding multiplicative seminorm $\zeta_{a,r}(f)$ mapping $f = \sum_n c_n (x-a)^n$ to 
\[\zeta_{a,r}(f) = \max_n \{|c_n| r^n\}.\]
It can be described as a supremum norm over the disc in $K$ around $a$ with radius $r$. The point $\zeta_{a,r}$ in $(\mathbb{A}^1)^\an$ is called of type two, if the radius $r$ lies in the value group $|K^\ast|$. Otherwise it is called of type three. %Points of type four are realized as limits of series of points of type three and two. They do only appear over fields that are not spherically complete. 
Note that $\zeta_{a,r} = \zeta_{b,s}$ if and only if $|a-b| \leq  r = s$. 

The Berkovich projective line is an $\mathbb{R}$-tree  in the sense of \cite{baru}, Appendix B. It is uniquely arcwise connected. We denote the unique path between two elements $x, y$ of $(\mathbb{P}^1)^\an$ by $[x,y]$. For $a,b \in \mathbb{A}^1(K)$ and $r<s$ the path $[\zeta_{a,r}, \zeta_{a,s}]$ consists of all points $\zeta_{a, t}$ such that $t\in [r,s]$.  To determine $[\zeta_{a,r}, \zeta_ {b,s}]$
in general, put $R = \max\{r,s,|a-b|\}$ and note that $\zeta_{a,R} = \zeta_{b,R}$. Then $[\zeta_{a,r}, \zeta_ {b,s}] = [\zeta_{a,r}, \zeta_{a,R}] \cup [\zeta_{b,R}, \zeta_{b,s}]$. 

The complement of the subset of points of type 1 in $(\mathbb{A}^1)^\an$ can be endowed with a path distance metric, see \cite{baru}, section 2.7. For points $x = \zeta_{a,r}$ and $y = \zeta_{b,s}$ as above, it is equal to 
\[\rho(x,y) = 2 \log R - \log r - \log s.\]
Set $\gaussBi=\zeta_{\zip,\rip}$ and $\gaussCi=\zeta_{\zim,\rim}$. 
Note that $\gamma_i(\gaussCi) = \gaussBi$. By the explicit description above,  is easy to see that the path $[\gaussCi,  \gaussBi[$ is contained in the fundamental domain $\mF$.

The Berkovich curve $X^\an$ contains the skeleton $\skeleton(X^\an)$ of the semistable model of $X$ constructed by Mumford \cite{mu}, Theorem 3.3. The graph $\skeleton(X^\an)$  is equal to the
image under the natural map $\pi:\Omega\to X^\an$ of the convex hull of all the points $\gaussBi, \gaussCi$ in $(\PP^1)^\an$.  Note that this convex hull is contained in $(\Omega)^\an$. 
The skeleton $\skeleton(X^\an)$ can be endowed with a natural metric induced by the path distance metric on the projective line.  It contains the $g$ cycles $\cycle_i = \pi([\gaussCi, \gaussBi])$.

\Prop{length}{~
\begin{enumerate}[i)]
 \item The cycle $\cycle_i$ has length $-\log |q_{ii}|$. 
 \item The intersection $\cycle_i \cap \cycle_j$ has length $|\log|q_{ij}||$. 
\end{enumerate}
}

\Prf{
i) After exchanging $\gamma_i$ with $\gamma_i^{-1}$ if necessary, we may assume that $|\zip| < |\zim|$. The length of $\cycle_i$ is equal to the path distance $\rho(\gaussCi, \gaussBi)$. Since $\rip,\ \rim < |\zim| = |\zim - \zip|$, we find 
\[\rho(\gaussCi, \gaussBi) = 2 \log |\zim| - \log \rip - \log\rim.\]
On the other hand, $|q_{ii}| = |u_i(z) / u_i(\gamma_iz)|$ for any point $z \in \Omega$. 
Let us choose a rational point $z \in \partial \discB_i$. 
Then $\gamma_iz$ lies in $\partial \discC_i$. Therefore we can apply Proposition \ref{prop:ufunction} to $z$ and $\gamma_iz$ to deduce
\[-\log |q_{ii}| = -\log \frac{|z-\zip||\gamma_i z - \zim|}{|z- \zim| |\gamma_i z - \zip|}.\]
Since by the choice of $z$ we have $|z - \zip|=\rip$ and $|\gamma_i z - \zim| = \rim$ and $ |z - \zim| = | \zim|$ as well as $|\gamma_i z - \zip| = |\gamma_i z | = |\zim|$, we deduce
that
\[  -\log |q_{ii}| = 2 \log |\zim| -  \log \rip -  \log \rim .\]

ii) After exchanging $\gamma_i$ and $\gamma_i^{-1}$ if necessary, we may again assume that $|\zip| < |\zim|$. Similarly, we may assume that  $|\zjp| < |\zjm|$. After exchanging $i$ and $j$ if necessary we also have $|\zip| \leq |\zjp|$. 

Choose $z \in \partial \discB_j$. Then $\gamma_jz \in \partial \discC_j$. 
Hence we can apply  Proposition \ref{prop:ufunction} to $z$ and $\gamma_jz$ to deduce
\[-\log |q_{ij}| = - \log \left|\frac{u_i(z)}{u_i(\gamma_j z)}\right| =  -\log \frac{|z-\zip||\gamma_j z - \zim|}{|z- \zim| |\gamma_j z - \zip|}.\]
By the choice of $z$ we have $|z - \zip| = |\zjp- \zip|$ and $ |z - \zim| = | \zjp -\zim|$ as well as $|\gamma_jz - \zip| = |\zjm - \zip |$ and $|\gamma_j z - \zim|= |\zjm - \zim|$. Hence
\begin{equation}\label{eq:length} -\log |q_{ij}|  = \log \frac {| \zjp -\zim| \, |\zjm - \zip |}{|\zjp- \zip|\, | \zjm -\zim|}.\end{equation}
Recall that $[\gaussBi, \gaussCi] = [\gaussBi, \zeta_{\zip, |\zim|}] \cup [\zeta_{\zim, |\zim|}, \gaussCi]$. The first part $[\gaussBi, \zeta_{\zip, |\zim|}]$ intersects $[\zeta_j^+, \zeta_j^-]$ in all points $\zeta_{\zip, r}$ satisfying $|\zip - \zjp| \leq r \leq  \min\{|\zim|, |\zjm|\}$.
The second part $[\zeta_{\zim,|\zim|}, \gaussCi]$ intersects $[\zeta_j^+, \zeta_j^-]$ in all points $\zeta_{\zim, r}$ such that either $|\zjm - \zim| \leq r \leq \min \{|\zim|, |\zjm|\}$
 or $|\zjp - \zim|  \leq r \leq \min\{|\zim|,|\zjm|\}$. 

Now we can check our claim case by case. 

{\bf Case 1:} $|\zim| < |\zjp|$. This implies $|\zip| < |\zim| < |\zjp| < |\zjm|$. In this case $\cycle_i \cap \cycle_j$ is empty and by formula \eqref{eq:length} we have $- \log|q_{ij}|  = 0$. 

{\bf Case 2:} $|\zim| = |\zjp|$. This implies $|\zip| < |\zim| = | \zjp| < |\zjm|$. Here, $\cycle_i \cap \cycle_j$ is the image of $[\zeta_{\zim,|\zjp - \zim|}, \zeta_{\zim,  |\zim|}]$, which has length $\log ( | \zim| | \zjp- \zim|\inv)$, which coincides with $ \log|q_{ij}| = \left|-\log|q_{ij}|\right|$ by formula \eqref{eq:length}.

{\bf Case 3:}  $|\zim| > |\zjp|$ and $|\zim| \neq |\zjm|$.   
Put $m = \min \{|\zim |, |\zjm|\}$. Then $\cycle_i \cap \cycle_j$ is the image of $[\zeta_{\zip, |\zip - \zjp| }, \zeta_{\zip, m}]$, which has length $\log ( m |\zip - \zjp|\inv)$. This coincides with $- \log|q_{ij}| $ by formula \eqref{eq:length}.

{\bf Case 4:} $|\zim| > |\zjp|$ and $|\zim| = |\zjm|$. In this case, $\cycle_i \cap \cycle_j$ is the image of $[\zeta_{ \zip, |\zip - \zjp| }, \zeta_{\zip, |\zim|}] \cup [\zeta_{\zim, |\zim|}, \zeta_{\zim, |\zim - \zjm|}]$, which has length
\[\log \frac{|\zim| |\zim|}{|\zip - \zjp| |\zim - \zjm|}  = - \log |q_{ij}|.\]

}

Note that part ii) of the preceding statement follows also from \cite{vp}, proof of Theorem 6.4 (2).

We can use the functions $u_i$ to give an explicit description of the embedding of $\curve$ into its Jacobian.

Fix $a \in \Omega(K)$ and consider 
 the analytic morphism
\[u: \Omega \longrightarrow  (\mathbb{G}_m^g)^\an, \quad z \mapsto \left(\frac{u_1(z)}{u_1(a)}, \ldots, \frac{u_g(z)}{u_g(a)} \right).\]

Let $\Lambda$ be the multiplicative lattice in $(\mathbb{G}_m^g)^\an$ generated by all $\lambda_i = (q_{i1}, \ldots, q_{ig})$. 

Then the quotient $\quotient{(\mathbb{G}_m^g)^\an}{\Lambda}$ is isomorphic to the analytification $J^\an$ of the Jacobian $J$ of $X$, and the following uniformization diagram is commutative, where $j$ is the embedding of $X$ into its Jacobian associated to the point $\pi(a) \in X(K)$.

\begin{center}
 \begin{tikzpicture}
  \matrix (m) [matrix of math nodes,row sep=3em,column sep=4em,minimum width=2em]
  {
     \Omega^\an & (\GG_m^g)^\an \\
     X^\an & J^\an \\};
  \path[-stealth,right hook->]
    (m-2-1) edge node [above] {$j^\an$} (m-2-2)
    (m-1-1) edge node [above] {$u$} (m-1-2);
  \path[-stealth]
    (m-1-1) edge (m-2-1)
    (m-1-2) edge (m-2-2);  
\end{tikzpicture}
\end{center}

There exists a birational morphism $j^{(g)}:X^{(g)}\rightarrow J$ where $X^{(g)}$ denotes the $g$-th symmetric power of the curve sending $(z_1,\dots,z_g)$ to $j(z_1)\cdots j(z_g)$. 

There is a natural tropicalization map 

\[\trop:  (\mathbb{G}_m^g)^\an \longrightarrow \mathbb{R}^g,  \]
which maps a multiplicative seminorm $x$ on the Laurent polynomial ring $K[T_1^{\pm 1}, \ldots, T_g^{\pm 1}]$ to $(\log|T_1(x)|, \ldots, \log|T_g(x)|)$, where $|T_i(x)|$ denotes the value of $x$ on $T_i$. 

It induces a tropicalization map $\trop: J^\an \longrightarrow \quotient{\mathbb{R}^g}{\log|\Lambda|}$.
We will use this tropicalization map on the Jacobian to construct faithful tropicalizations.

\section{Mumford curves of genus two}
We study the situation of Mumford curves of genus $g=2$. In this case, the skeleton $\skeleton(X^\an)$ contains two cycles $\cycle_1$ and $\cycle_2$. 
These cycles are either disjoint, meet in a point or meet in an edge (figure \ref{fig:skeletons}). We know from the proof of Proposition \ref{prop:length} that in the first two cases $|q_{12}| = |q_{21}| = 1$, and that in the third case the length of the shared edge is $|\log|q_{12}||$.

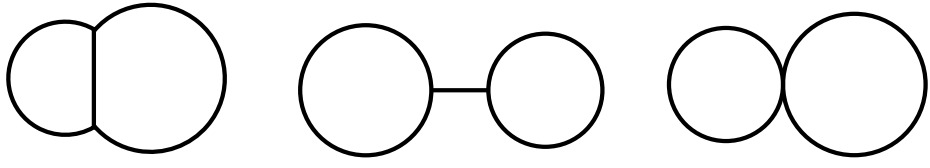
\fig{skeletons}{The three possible types of skeletons in genus two, from the left to the right: Skeleton with shared edge, connecting edge and connecting point.}{
 \begin{tikzpicture}[scale=.75]
  \draw[ultra thick] (0,0) arc (60:300:1) -- ++(0,1.7) arc (320:40:-1.3);
 \end{tikzpicture}
 \qquad
 \begin{tikzpicture}[scale=.75]
   \draw[ultra thick] (-1.65,0) circle (1.15);
  \draw[ultra thick] (+1.5,0) circle (1);
  \draw[ultra thick] (-.5,0) -- (.5,0);
 \end{tikzpicture}
 \qquad
  \begin{tikzpicture}[scale=.75]
  \draw[ultra thick] (-1,0) circle (1);
  \draw[ultra thick] (+1.25,0) circle (1.25);
 \end{tikzpicture}
}

\Lemma{sortieren}{
For a Schottky group $\Gamma$ of rank two there exists an element $\gamma \in \PGL(2,K)$ such that $\gamma \Gamma \gamma^{-1}$ has a good fundamental domain satisfying 
$|\zeinsp|~<~|\zeinsm|~<~|\zweim|$ and $|\zeinsp|~<~|\zweip|~<~|\zweim|$. If the skeleton has a shared edge, we may also assume $|\zweip|<|\zeinsm|$. 
}
\Prf{Let $\gamma_1, \gamma_2 $ be generators of $\Gamma$ and $\discB_1,\discB_2,\discC_1,\discC_2$ discs in $\PP^1$ providing a good fundamental domain. For $i= 1,2$ let  $\zip$ be the fixed point of $\Gamma$ in $\discB_i$, and let $\zim$ be the fixed point of $\Gamma$ in $\discC_i$. Recall that all elements in $\discB_i$ have absolute value $|\zip|$, all elements in $\discC_i$ have absolute value $|\zim|$, and the discs $\discB_1,\discB_2,\discC_1, \discC_2$ are pairwise disjoint. We have 
$$|\zeinsm-\zeinsp|\leq\max\{|\zeinsm-\zweim|,|\zweim-\zeinsp|\}.$$
If $|\zeinsm-\zweim|>|\zweim-\zeinsp|$ 
we exchange the generators $\gamma_1$ and $\gamma_1^{-1}$, and hence the pair of fixed points. Therefore we may assume that 
 $|\zeinsm-\zweim|\leq |\zweim-\zeinsp|$, which implies $|\zeinsm - \zeinsp| \leq |\zweim - \zeinsp|$.  Similarly, after exchanging $\gamma_2$ and $\gamma_2^-$ if necessary, we may assume that  $|\zeinsp-\zweip|\leq|\zweim-\zeinsp|$, which implies that $|\zweip-\zweim| \leq |\zweim - \zeinsp|$.
Since
$$|\zeinsp-\zweip|\leq\max\{|\zeinsp-\zweim|,|\zweim-\zweip|\}$$ we have $|\zeinsp-\zweip|\leq|\zweim-\zeinsp|$ and similarly $|\zeinsm-\zweim|\leq |\zweim-\zeinsp|$. 

Thus $$|\zweim-\zeinsp|\geq \max\{|\zeinsm-\zeinsp|,|\zweim-\zweip|,|\zweim-\zeinsm|,|\zweip-\zeinsp|\}.$$ 

We find a point $p \in K$ lying outside the closed disc $\discB_1$ such that $|\zeinsp-p|<\min \{|\zeinsp-\zeinsm|,|\zeinsp-\zweip|,|\zweim -\zeinsp|\}$ 
and there exists a point $q \in K$ lying outside the closed disc $\discC_2$ satisfying $|\zweim-q|<\min \{|\zweim-\zeinsp|,|\zweim-\zeinsm|,|\zweim-\zweip|\}$. Such points do exist since the discs are pairwise disjoint. 
This implies that 
$$|\zeinsm-p|=|\zeinsm - \zeinsp|, |\zweip-p| =|\zweip - \zeinsp| \mbox{ and } |\zweim - p| =  |\zweim-\zeinsp|,$$
as well as 
$$|\zeinsp-q|=|\zeinsp-\zweim|, |\zeinsm - q| = |\zeinsm-\zweim| \mbox{ and } |\zweip - q| = |\zweip - \zweim|.$$
Hence $p$ and $q$ lie in the fundamental domain $\mathcal{F}$. 

We conjugate $\Gamma$ by the element $\gamma(z) = \frac{z-p}{z-q}$ in $\PGL(2,K)$. Then we obtain
\begin{align*}
 |\gamma(\zeinsp)|&=\left|\dfrac{\zeinsp-p}{\zeinsp-q}\right|
 <\left|\dfrac{\zeinsm-p}{\zeinsm-q}\right|=|\gamma(\zeinsm)|
\end{align*}
since $|\zeinsp-p|<|\zeinsm-p|$ and $|\zeinsm-q|=|\zweim-\zeinsm|\leq|\zweim-\zeinsp|=|\zeinsp-q|$.

Moreover, we have 
\begin{align*}
 |\gamma(\zeinsm)|&=\left|\dfrac{\zeinsm-p}{\zeinsm-q}\right| 
 <\left|\dfrac{\zweim-p}{\zweim-q}\right|=|\gamma(\zweim)|
\end{align*}
since $|\zweim-q|<|\zeinsm-q|$ and $|\zeinsm-p|=|\zeinsm-\zeinsp|\leq|\zweim-\zeinsp|=|\zweim-p|$.

This shows $|\gamma(\zeinsp)|<|\gamma(\zeinsm)|<|\gamma(\zweim)|$. 
We also find
\begin{align*}
 |\gamma(\zeinsp)|&=\left|\dfrac{\zeinsp-p}{\zeinsp-q}\right|
 <\left|\dfrac{\zweip-p}{\zweip-q}\right|=|\gamma(\zweip)|\end{align*}
since $|\zeinsp-p|<|\zweip-p|$ and $|\zeinsp-q|=|\zweim-\zeinsp|\geq|\zweip-\zweim|=|\zweip-q|$. Similarly,
\begin{align*}
 |\gamma(\zweip)|&=\left|\dfrac{\zweip-p}{\zweip-q}\right|
 <\left|\dfrac{\zweim-p}{\zweim-q}\right|=|\gamma(\zweim)|
\end{align*}
since $|\zweip-p|=|\zweip-\zeinsp|\leq|\zweim-\zeinsp|=|\zweim-p|$ and $|\zweip-q| = |\zweip - \zweim|>|\zweim-q|$.
Thus we also get $|\gamma(\zeinsp)|<|\gamma(\zweip)|<|\gamma(\zweim)|$.

Hence $\gamma \Gamma \gamma^{-1}$ has a good fundamental domain satisfying the first claim of the lemma.

In order to prove the second claim, we consider the case that the skeleton has a shared edge. By the first part of the proof  we may assume that $|\zeinsp|< |\zeinsm|<|\zweim|$   and  $|\zeinsp|<|\zweip| < |\zweim|$ hold. Recall from the proof of Proposition \ref{prop:length} that the inequality   $|\zweip|>|\zeinsm|$ implies that the cycles are  disjoint. Hence we have $|\zweip| \leq |\zeinsm|$.  If the inequality is strict, our claim follows. Therefore we may assume that $|\zeinsm|=|\zweip|$.
Note that in the case of a shared edge we have $|\zeinsm - \zweip| < |\zeinsm|  = |\zweip|$. This follows from the discussion of Case 2 in  the proof of Proposition \ref{prop:length}.

There exists a point $p\in K$ outside the disc $\discC_1$ such that 
\[|\zeinsm-p|<\min\{|\zweim - \zeinsm|, |\zeinsp - \zeinsm|,|\zweip-\zeinsm|\}\] and  a point $q\in K$ outside $\discC_2$ such that  
\[|\zweim-q| < \min \{|\zweim - \zeinsm|, |\zweim - \zweip|, |\zweim - \zeinsp|\}.\]
A similar argument as in the proof of the first claim shows that $p$ and $q$ lie in the fundamental domain $\mathcal{F}$. Let $\gamma(z)=\frac{z-p}{z-q}$. We find
\begin{align*}
 |\gamma(\zeinsm)|&=\left|\dfrac{\zeinsm-p}{\zeinsm-q}\right|=\dfrac{|\zeinsm-p|}{|\zweim|} 
 <\dfrac{|\zweip-p|}{|\zweim|}=\left|\dfrac{\zweip-p}{\zweip-q}\right|=|\gamma(\zweip)|
\end{align*}
since $|\zeinsm-p|<|\zweip - \zeinsm| = |\zweip-p|$ and $|\zeinsm-q|=|\zeinsm - \zweim| = |\zweim|$ as well as $|\zweip-q|= |\zweip - \zweim| = |\zweim|$. 
Similarly, we find
\begin{align*}
 |\gamma(\zweip)|&=\left|\dfrac{\zweip-p}{\zweip-q}\right|=\dfrac{|\zweip-\zeinsm|}{|\zweim|} 
 <\dfrac{|\zeinsp-\zeinsm|}{|\zeinsp-q|}=\left|\dfrac{\zeinsp-p}{\zeinsp-q}\right|=|\gamma(\zeinsp)|
\end{align*}
since $|\zweip - p| = |\zweip-\zeinsm|<|\zeinsm| = |\zeinsp-\zeinsm| = |\zeinsp - p|$ and $|\zweip-q|=|\zweip - \zweim|=|\zweim| = |\zweim - \zeinsp| = |\zeinsp - q|$.
Moreover, we have 
\begin{align*}
 |\gamma(\zeinsp)|&=\left|\dfrac{\zeinsp-p}{\zeinsp-q}\right|=\dfrac{|\zeinsm|}{|\zweim|} 
 <\dfrac{|\zweim|}{|\zweim-q|}=\left|\dfrac{\zweim-p}{\zweim-q}\right|=|\gamma(\zweim)|
 \end{align*}
since $|\zeinsp - p| = |\zeinsp - \zeinsm| = |\zeinsm|<|\zweim|$ and $|\zweim-q|<|\zweim|$.

Thus $|\gamma(\zeinsm)|<|\gamma(\zweip)|<|\gamma(\zeinsp)|<|\gamma(\zweim)|$.
Now exchange $\gamma_1$ and $\gamma_1^{-1}$ which exchanges $\zeinsm$ and $\zeinsp$. This concludes the proof. 
}

Let $\curve$ be a Mumford curve of genus $2$ over $K$ obtained by a Schottky uniformization with a Schottky group $\Gamma$ admitting a good fundamental domain. Since conjugation of the Schottky group by an element in $\PGL(2,K)$ amounts to a $K$-rational automorphism of the curve, we may assume by Lemma \ref{lemma:sortieren} that the set of generators of $\Gamma$ satisfy the inequalities $|\zeinsp| < |\zeinsm| < |\zweim|$ and $|\zeinsp| < |\zweip| < |\zweim|$ with the additional inequality $|\zweip| < |\zeinsm|$ in the case of a shared edge.

Every $K$-rational point in $\curve$ (and hence every $K$-rational point of $\Omega$) induces an embedding $j: \curve \hookrightarrow J$ into the Jacobian. Recall from the last section that the uniformization of $J^\an$ by a two-dimensional torus gives rise to a tropicalization map $\trop: J^\an \rightarrow \quotient{\RR^2}{\log |\Lambda|}$. 
We consider the composition
\[ \mu: \skeleton(\curve^\an) \hookrightarrow  \curve^\an \stackrel{j^\an}{\longhookrightarrow} J^\an \stackrel{\trop}{\longrightarrow} \quotient{\RR^2}{\log|\Lambda|}\] 
Note that $\mu$ is in general not injective.

\Thm{tropicalimage}{~
Denote by $e_1, e_2$ the canonical basis of $\mathbb{R}^2$. 
There exists a point $a \in \Omega(K)$ such that the map $\mu$ defined with the embedding $j: X \hookrightarrow J$ associated to $a$ has
 the following image in a fundamental parallelepiped in $\quotient{\mathbb{R}^2}{\log|\Lambda|}$:

\begin{enumerate}
\item  In the case of disjoint cycles or cycles meeting in a point we have \[\im(\mu)= [0,-\log|q_{11}|[ \, e_1 \cup [0, - \log |q_{22}|[ \, e_2.\]
 \item In the case of cycles sharing an edge, we put $l = |\log|q_{12}|| = |\log|q_{21}||$ and $v = (l,l)$. Then we have 
 \[ \im(\mu)=[0,1[\, v\cup v+[0,\log|q_{11}| -l[\, e_1\cup v+[0,\log|q_{22}|-l[ \, e_2.\]
 
\item In the case of cycles sharing an edge, we have furthermore $\mu(\cycle_1\setminus \cycle_2)\subset v+[0,\log|q_{11}| - l [\, e_1$ and  $\mu(\cycle_2\setminus \cycle_1)\subset v+[0,\log|q_{22}| -l[ \, e_2$. 

 \item $\mu$ maps the two cycles $\cycle_1,\cycle_2$ in $\skeleton(\curve^\an)$ isometrically to their images $\mu(\cycle_1)$ resp. $\mu(\cycle_2)$ where we endow the fundamental parallelepiped with the metric induced from the maximum metric on $\RR^2$. 
\end{enumerate}
}

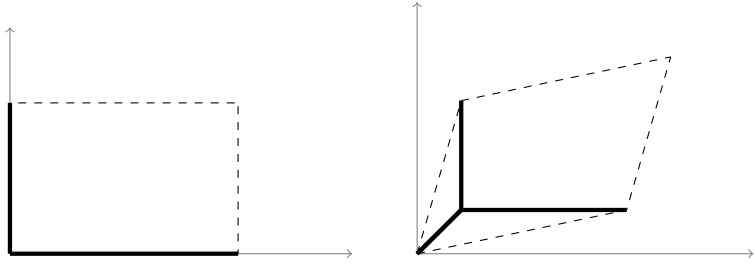
\fig{tropicalimage}{The fundamental parallelepiped and the image of $\mu$. Depicted on the left is the case of a connecting edge, on the right the case of a shared edge.}{

\begin{tikzpicture}[scale=2, m/.style={ultra thick, black}]
\def\lam{0} %Länge der shared edge
\def\c{1.5} %Länge Zykel 1
\def\d{1} %Länge Zykel 2

\coordinate (v) at (\lam,\lam);
\coordinate (C1) at (\c,0);
\coordinate (C2) at (0,\d);

\coordinate (Z) at (0,0);

%Koordinatenachsen
\draw[gray,->] (Z) --++(0,1.5*\d+2*\lam);
\draw[gray,->] (Z) --++(1.5*\c+2*\lam,0);

%Fundamentalbereich
 \draw[dashed] (Z) --++(\c+\lam,\lam) --++(\lam,\lam+\d);
 \draw[dashed] (Z) --++(\lam,\d+\lam) --++(\c+\lam,\lam);
 %Bild von \mu
 \draw[m] (Z)--++(v)--++(C1);
 \draw[m] (Z)--++(v)--++(C2);
 \end{tikzpicture}
\qquad
\begin{tikzpicture}[scale=1.45, m/.style={ultra thick, black}]
\def\lam{.4} %Länge der shared edge
\def\c{1.5} %Länge Zykel 1
\def\d{1} %Länge Zykel 2

\coordinate (v) at (\lam,\lam);
\coordinate (C1) at (\c,0);
\coordinate (C2) at (0,\d);

\coordinate (Z) at (0,0);

%Koordinatenachsen
\draw[gray,->] (Z) --++(0,1.5*\d+2*\lam);
\draw[gray,->] (Z) --++(1.5*\c+2*\lam,0);

%Fundamentalbereich
 \draw[dashed] (Z) --++(\c+\lam,\lam) --++(\lam,\lam+\d);
 \draw[dashed] (Z) --++(\lam,\d+\lam) --++(\c+\lam,\lam);
 %Bild von \mu
 \draw[m] (Z)--++(v)--++(C1);
 \draw[m] (Z)--++(v)--++(C2);
 \end{tikzpicture}
}

\Prf{Recall that by Proposition \ref{prop:ufunction} $|u_i(z)| = \frac{|z-\zip|}{|z-\zim|}$ on the fundamental domain. 
Hence we find for all $x = \zeta_{\zeinsp, r}$ with $\reinsp \leq r \leq |\zeinsm|$:

\[(|u_1(x)|, |u_2(x)| )
 =\left(\frac{r}{|\zeinsm|}, \frac{\max\{r, |\zeinsp - \zweip|\}}{|\zweim|}\right) \]

For all $x = \zeta_{\zeinsm, r}$ with $\reinsm \leq r \leq |\zeinsm|$ we have
\[
(|u_1(x)|, |u_2(x)| )
 = \left(\frac{|\zeinsm|}{ r }, \frac{\max\{r,|\zeinsm - \zweip|\}}{\max\{r,|\zeinsm - \zweim|\}}\right) \]

Let us first discuss the case where the cycles are disjoint or share a point. By Proposition \ref{prop:length} this correponds to the case $|\zeinsp| < |\zeinsm| < |\zweip| < |\zweim|$ or
to the case $| \zeinsp| < |\zeinsm| = |\zweip| < |\zweim|$ with $|\zweip - \zeinsm| = |\zeinsm|$. 
Choose a point $a \in \Omega(K)$ such that $|a| =|\zeinsm| =  |\zeinsm -a|$ and $|\zweip| =  | \zweip -a|$. 
Then
$(|u_1(a)|, |u_2(a)|)  = ( 1, |\zweip| /|\zweim|)$. If we use $\pi(a)$ to define the embedding $j$ of $\curve$ into its Jacobian, the associated map $\mu$ is given by
$\mu(x) = \left(\log| u_1(x)|/| u_1(a)|, \log |u_2(x)|/ |u_2(a)|\right) \mod \log |\Lambda|$ on the fundamental domain. Hence
the cycle $\cycle_1$ is mapped to
\[ \left\{(\log r - \log|\zeinsm|,0 ) : r \in [\reinsp, |\zeinsm|]\right\} \cup \left\{ (\log|\zeinsm | - \log r, 0): r \in [\reinsm, |\zeinsm|]\right\}.\]
Recall that in this case $\log|q_{12}| = 0$, so that we can add $(-\log|q_{11}|, 0 ) = (2 \log |\zeinsm| - \log \reinsp-\log\reinsm,0) \in \log |\Lambda|$ to the first subset. Therefore the image of $\cycle_1$ is $\left[0,- \log|q_{11}|\right[ \, e_1$.  The second cycle $\cycle_2$ can be treated in the same way. This proves claim (1). 

Now let us discuss the cases of cycles sharing an edge. In this case we have $|\zeinsp| < |\zweip| < |\zeinsm| < |\zweim|$. Choose a point $a \in \Omega(K)$ such that $|a| = |\zweip| = |\zweip - a|$. Then $(|u_1(a)|, |u_2(a)|) = (|\zweip| / |\zeinsm |,  |\zweip| / |\zweim|)$. If we use $\pi(a)$ to define the embedding $j$ of $\curve$ into its Jacobian, the associated map $\mu$ is given by
$\mu(x) = (\log| u_1(x)|/| u_1(a)|, \log |u_2(x)|/ |u_2(a)|) \mod \log |\Lambda|$ on the fundamental domain.
Recall from the proof of Proposition \ref{prop:length} that the intersection $\cycle_1 \cap \cycle_2$ corresponds to the line segment $[\zeta_{\zeinsp, |\zweip|}, \zeta_{\zeinsp, |\zeinsm|}]$. Hence $\cycle_1 \cap \cycle_2$ is  mapped to the line segment between $(0,0)$ and $(\log|\zeinsm| - \log|\zweip|, \log|\zeinsm| - \log|\zweip|)$ with slope one. Formula \eqref{eq:length} in the proof of Proposition \ref{prop:length} implies that $- \log|q_{12}| = \log |\zeinsm| - \log |\zweip|$. Therefore the shared edge is mapped to $[0,1] \, v$. Let us now look at $\cycle_1 \backslash \cycle_2$. It corresponds to the union of the paths $[\zeta_{\zeinsp, \reinsp}, \zeta_{\zeinsp, |\zweip|}]$ and $[\zeta_{\zeinsm, \reinsm}, \zeta_{\zeinsm, |\zeinsm|}]$. 
On the first path $[\zeta_{\zeinsp, \reinsp}, \zeta_{\zeinsp, |\zweip|}]$ the map $\mu$ is given by \[\zeta_{\zeinsp, r} \mapsto \left(\log \frac{ r}{|\zweip|}, 0 \right).\]
Adding $(\log|q_{11}|, \log|q_{12}| ) = (-2 \log |\zeinsm| + \log \reinsp +\log \reinsm, -\log|\zeinsm| + \log|\zweip|) \in \log |\Lambda|$ we find that this is the line segment between
$(\log\frac{|\zeinsm|^2}{\reinsm |\zweip|}, l)$ and the lattice point $(\log|q_{11}|, l)$.
On the second path $[\zeta_{\zeinsm, \reinsm}, \zeta_{\zeinsm, |\zeinsm|}]$ the map $\mu$ is given by 
\[ \zeta_{\zeinsm, r} \mapsto \left(\log\frac{|\zeinsm|^2}{r |\zweip|}, \log \frac{|\zeinsm|}{|\zweip|} \right).\]
Its image is the line segment between $v$ and $(\log\frac{|\zeinsm|^2}{\reinsm |\zweip|}, l)$. The cycle $\cycle_2$ can be treated similarly. 
Hence claims (2) and (3) hold.
To show claim (4), write $\mu(x)=(\mu_1(x),\mu_2(x))$. Recall that on the first cycle $\mu_1$ is of the form $\log r$ plus constant. In the case of a connecting edge $\mu_2$ is equal to zero, which implies that $\mu$ is an isometry on this cycle. In the shared edge  case there is also a change in $\mu_2$, but $\mu_1$ will increase in the same way, so the maximum will be achieved in $\mu_1$. A similar argument holds for the second cycle.

}

\Lemma{twoadd}{
Again let $v=(-\log|q_{12}|,-\log|q_{12}|)$. If $x=2v+\alpha e_1+\beta e_2$ with $0<\alpha<-\log|q_{11}|+\log|q_{12}|$ and $0<\beta<-\log|q_{22}|+\log|q_{12}|$ then there are exactly two points $S$ and $T$ in the skeleton such that $\mu(S)+\mu(T)=x$, namely $S=\mu^{-1}(v+\alpha e_1)$ and $T=\mu^{-1}(v+\beta e_2)$. 
}
\Prf{
If $\mu(S)=v+\alpha e_1$ and $\mu(T)=v+\beta e_2$ then clearly $\mu(S)+\mu(T)=x$. So let us suppose there exist two other points $U$ and $V$ in the skeleton such that $\mu(U) + \mu (V)=x$. We have to handle several cases:
\begin{enumerate}
 \item $\mu(U)=v+\alpha' e_1$ and $\mu(V)=v+\beta' e_2$ with $0<\alpha'<-\log|q_{11}|+\log|q_{12}|$ and $0<\beta'<-\log|q_{22}|+\log|q_{12}|$. Then $2v+\alpha e_1+\beta e_2=2v+\alpha' e_1+\beta' e_2$ so $v+\alpha e_1+\beta e_2=v+\alpha' e_1+\beta' e_2$. Since these points lie in the fundamental domain of the tropicalization of the Jacobian we obtain $\alpha'=\alpha$ and $\beta=\beta'$, thus $U$ and $V$ coincide with $S$ and $T$. 
 \item $\mu(U)=v+\alpha' e_1$ and $\mu(V)=v+\alpha'' e_1$ with $0<\alpha',\alpha''<-\log|q_{11}|+\log|q_{12}|$. Then 
 $\mu(U) + \mu(V)$ lies inside the image of the interior of the first cycle, but $x$ does not. Thus such points $U$ and $V$ cannot exist. The same argument holds if both points lie on the other cycle. 
 \item $\mu(U)=v+\alpha' e_1$ and $\mu(V)=\gamma v$ with $0<\alpha'<-\log|q_{11}|+\log|q_{12}|$ and $\gamma\in [0,1]$. Then
 $2v+\alpha e_1+\beta e_2=(1+\gamma)v+\alpha' e_1$ thus $v+\alpha e_1+\beta e_2=\gamma v+\alpha' e_1
 $.
 To satisfy this equality $\gamma v+\alpha' e_1$ cannot lie in the fundamental domain of the tropicalization of the Jacobian. 
 Then $(1+\gamma)v+\alpha'e_1-\log|q_{22}|e_2$ lies in the fundamental domain, but this point is not equal to $x=2v+\alpha e_1+\beta e_2$.
 \item $\mu(U)=\gamma v$ and $\mu(V)=\gamma' v$ with $\gamma,\gamma'\in [0,1]$. Then $2v+\alpha e_1+\beta e_2=(\gamma'+\gamma)v$ thus $v+\alpha e_1+\beta e_2=\delta v,\ \delta\in [-1,1] 
 $. This is not possible if $\delta\geq 0$. So suppose $\delta< 0$. 
We replace $\delta v$ by an equivalent point:
 \begin{align*}
 \delta v&=2v+(-\log|q_{11}|+\log|q_{12}|)e_1+(-\log|q_{22}|+\log|q_{12}|)e_2+\delta v \\
&=(1+\delta)v-\log|q_{11}|e_1-\log|q_{22}|e_2
\end{align*}
Thus $\delta v-\log|q_{11}|e_1-\log|q_{22}|e_2=\alpha e_1+\beta e_2$ and we obtain $\alpha=-\log|q_{11}|-\delta\log|q_{12}|$. This contradicts $\alpha<-\log|q_{11}|+\log|q_{12}|$. %

The same argument holds for $\beta$ thus the points $U$ and $V$ cannot exist.

\end{enumerate}
}

\section{Faithful tropicalization in genus two}

In this section we will always assume that the cycle lengths are at least twice as long as the length of the shared edge, i.e. $|q_{11}|<|q_{12}|^2$ and $|q_{22}|<|q_{12}|^2$.

\Lemma{points}{
Let $A=-\log|q_{11}|+2\log|q_{12}|$, $B=-\log|q_{22}|+2\log|q_{12}|$ and\\ $v~=~(-\log|q_{12}|,-\log|q_{12}|)$. %
Then there exist points \\ $P_1,P_2,P_3,P_4,S_1,S_2,S_3,T_1,T_2,T_3\in~X(K)$ such that
\begin{itemize}
\item $\mu\circ\tau (P_1)=v+\frac{1}{4}A e_1$
\item  $\mu\circ\tau (P_2)=v+\frac{1}{2}A e_1$
\item $\mu\circ\tau (P_3)=v+\frac{1}{4}B e_2$
\item $\mu\circ\tau (P_4)=v+\frac{1}{2}B e_2$
\item $\mu\circ\tau(S_1)=\mu\circ\tau(T_3)=v+(\frac{3}{4}B-\log|q_{12}|) e_2$
\item $\mu\circ\tau(T_1)=\mu\circ\tau(S_3)=v+(\frac{3}{4}A-\log|q_{12}|) e_1$
\item $\mu\circ\tau(S_2)=v-\frac{1}{2}\log|q_{22}|e_2$
\item $\mu\circ\tau(T_2)=v-\frac{1}{2}\log|q_{11}|e_1$ 
\item $\dfrac{j(P_1)j(P_2)j(P_3)}{j(S_1)j(S_2)j(S_3)},\dfrac{j(P_1)j(P_3)j(P_4)}{j(T_1)j(T_2)j(T_3)}$ are trivial in $J$.
\end{itemize}
}
\Prf{

All points on the skeleton stated in the lemma are of type two thus there exist $K$-rational points $P_1,P_2,P_3,P_4,S_1$ and $T_1$ as postulated. Since the map $j^{(2)}$ from $X^{(2)}$ to $J$ is surjective there are points $S_2,S_3$ and $T_2,T_3$ in $X(K)$ with
$$j(S_2)j(S_3)~=~\dfrac{j(P_1)j(P_2)j(P_3)}{j(S_1)} \text{ and } j(T_2)j(T_3)~=~\dfrac{j(P_1)j(P_3)j(P_4)}{j(T_1)}$$ thus $\dfrac{j(P_1)j(P_2)j(P_3)}{j(S_1)j(S_2)j(S_3)}$ and $\dfrac{j(P_1)j(P_3)j(P_4)}{j(T_1)j(T_2)j(T_3)}$ are trivial in $J$. To determine the position of the retractions of $S_2$ and $S_3$ on the skeleton we use that
\begin{align*}
\mu\circ\tau(S_2)+\mu\circ\tau(S_3)&=\mu\circ\tau(P_1)+\mu\circ\tau(P_2) + \mu\circ\tau(P_3) -\mu\circ\tau(S_1) \mod{\log|\Lambda|}.
\end{align*}
Thus
\begin{align*}
\mu\circ\tau(S_2)+\mu\circ\tau(S_3)&=2v+\frac{3}{4}A e_1-\left(\frac{1}{2}B+\log|q_{12}|\right)e_2  \mod{\log|\Lambda|} \\
&=2v+\frac{3}{4}A e_1+\frac{1}{2}\log|q_{22}|e_2 \mod{\log|\Lambda|} \\
&=3v+ \frac{3}{4}A e_1 + \left(-\log|q_{22}|+\log|q_{12}|+\frac{1}{2}\log|q_{22}|\right)e_2  \mod{\log|\Lambda|}
\end{align*}
where we replaced the point $\frac{1}{2}\log|q_{22}|e_2$ by an equivalent one. Splitting one of the vectors $v=-\log|q_{12}|e_1-\log|q_{12}|e_2$ yields
\begin{align*}
\mu\circ\tau(S_2)+\mu\circ\tau(S_3)&=2v + \left(\frac{3}{4}A-\log|q_{12}|\right)e_1 - \frac{1}{2}\log|q_{22}| e_2.
\end{align*}
Hence by Lemma \ref{lemma:twoadd}
$$\mu\circ\tau(S_2)=v-\frac{1}{2}\log|q_{22}| e_2 \text{ and }\mu\circ\tau(S_3)=v+\left(\frac{3}{4}A-\log|q_{12}|\right)e_1.
$$
$\mu\circ\tau(T_2)$ and $\mu\circ\tau(T_3)$ can be determined by a similar calculation.

}

\Thm{faithful_skeleton}{
Let $P_1,P_2,P_3,P_4\in \curve(K)$ be given as in the previous lemma and let $D=P_1+P_2+P_3+P_4$ be the corresponding effective divisor. 
Then there exist functions $f,g$ in the linear system $L(D)$ such that the rational morphism $(1,f,g)~:~\curve~\dashrightarrow ~\PP^2$ tropicalizes the skeleton faithfully.
}
\Prf{
By the previous lemma there exist points $S_1,S_2,S_3,T_1,T_2,T_3\in X(K)$ such that $S_1+S_2+S_3-P_1-P_2-P_3=\div(f)$ and  $~T_1~+~T_2+T_3-P_2-P_3-P_4=\div(g)$ are principal divisors. Both functions $f$ and $g$ are elements in the linear system $L(D)$. 
Since $\mu$ is an isometry on the cycles by theorem \ref{thm:tropicalimage} the points $P_1,P_2,P_3,P_4,S_1,S_2,S_3,T_1,T_2,T_3$ retract to the skeleton as shown in figure \ref{fig:extended_skeleton_f} where we label the edges on the first cycle by $a,b,c,d,e,\lambda$.
Note that $\lambda$ and $c$ are not present in the connecting edge case. We orient all edges so that $b,c,d,e,\lambda,a$ is the path $\tau(P_1)\rightarrow\tau(P_2)\rightarrow\tau(T_2)\rightarrow\tau(T_1)\rightarrow\tau(P_1)$.
Similarly we label the edges on the second cycle by $\alpha,\beta,\gamma,\delta$ and $\lambda$, where again $\lambda$ and $\gamma$ are not present in the case of a connecting edge. We orient the cycle so that $\beta,\gamma,\delta,\lambda,\alpha$ is the path $\tau(P_3)\rightarrow\tau(P_4)\rightarrow\tau(S_2)\rightarrow\tau(S_1)\rightarrow\tau(P_3)$.
In the connecting edge case the edge $\chi$ will be oriented from the first to the second cycle.

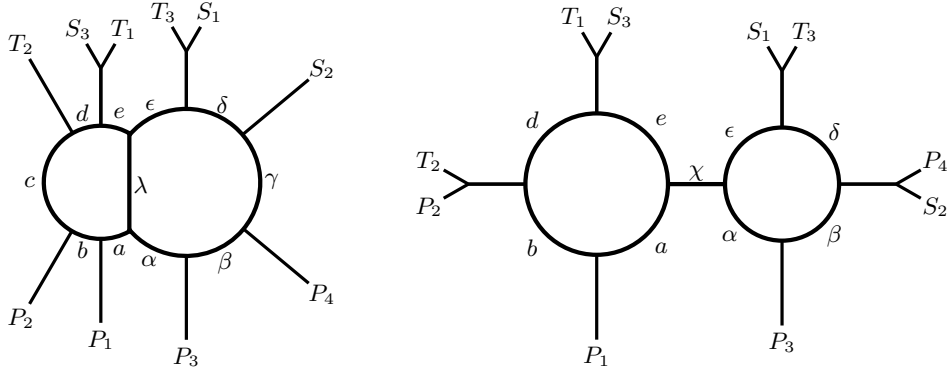
 \fig{extended_skeleton_f}{The extended skeletons.}{
 \begin{tikzpicture}[scale=.75]
%skeleton 
  \draw[ultra thick] (0,0) arc (60:300:1) -- ++(0,1.7) arc (320:40:-1.3);
 % \draw (1,-.85) circle (1.3);
%pts on left cycle
  \draw[very thick] (-.5,-.85) ++(270:1) -- +(270:1.5) +(270:1.8) node {\small $P_1$};
  \draw[very thick] (-.5,-.85) ++(240:1) -- +(240:1.5) +(240:1.8) node {\small $P_2$};
  \draw[very thick] (-.5,-.85) ++(120:1) -- +(120:1.5) +(120:1.8) node {\small $T_2$};
  \draw[very thick] (-.5,-.85) ++(90:1) -- ++(90:1) -- ++(120:0.5) ++(120:0.3) node {\small $S_3$};
  \draw[very thick] (-.5,-.85) ++(90:1) -- ++(90:1) -- ++(60:0.5) ++(60:0.3) node {\small $T_1$};
%pts on right cycle     
  \draw[very thick] (1,-.85) ++(270:1.3) -- +(270:1.5) +(270:1.8) node {\small $P_3$};
  \draw[very thick] (1,-.85) ++ (320:1.3) -- +(320:1.5) +(320:1.8) node {\small $P_4$};
  \draw[very thick] (1,-.85) ++(40:1.3) -- +(40:1.5) +(40:1.8) node {\small $S_2$};
  \draw[very thick] (1,-.85) ++(90:1.3) -- ++(90:1) -- ++(120:0.5) ++(120:0.3) node {\small $T_3$};
   \draw[very thick] (1,-.85) ++(90:1.3) -- ++(90:1) -- ++(60:0.5) ++(60:0.3) node {\small $S_1$};
%labels of edges
  \draw (.2,-.9) node {\small$\lambda$};
 \draw (1,-.85) ++ (245:1.55) node {\small$\alpha$};
  \draw (1,-.85) ++ (295:1.6) node {\small$\beta$};
  \draw (1,-.85) ++ (0:1.5) node {\small$\gamma$};
  \draw (1,-.85) ++ (65:1.5) node {\small$\delta$};
  \draw (1,-.85) ++ (115:1.5) node {\small$\epsilon$};
    \draw (-.5,-.85) ++ (285:1.25) node {\small$a$};
    \draw (-.5,-.85) ++ (255:1.25) node {\small$b$};
    \draw (-.5,-.85) ++ (180:1.25) node {\small$c$};
    \draw (-.5,-.85) ++ (105:1.25) node {\small$d$};
    \draw (-.5,-.85) ++ (75:1.25) node {\small $e$};
  \end{tikzpicture}
  \qquad
   \begin{tikzpicture}[scale=.75]
  \draw[ultra thick] (-1.75,0) circle (1.25);
  \draw[ultra thick] (+1.5,0) circle (1);
  \draw[ultra thick] (-.5,0) -- (.5,0);
  \draw[very thick] (-1.75,0) ++(270:1.25) -- +(270:1.5) +(270:1.8) node {\small $P_1$};
  \draw[very thick] (-2,0) ++(180:1) -- ++(180:1) -- ++(210:0.5) ++(210:0.3) node {\small $P_2$};
      \draw[very thick] (-2,0) ++(180:1) -- ++(180:1) -- ++(150:0.5) ++(150:0.3) node {\small $T_2$};
  \draw[very thick] (-1.75,0) ++(90:1.25) -- ++(90:1) -- ++(120:0.5) ++(120:0.3) node {\small $T_1$};
    \draw[very thick] (-1.75,0) ++(90:1.25) -- ++(90:1) -- ++(60:0.5) ++(60:0.3) node {\small $S_3$};
  \draw[very thick] (1.5,0) ++(270:1) -- +(270:1.5) +(270:1.8) node {\small $P_3$};
  \draw[very thick] (1.5,0) ++(90:1) -- ++(90:1) -- ++(120:0.5) ++(120:0.3) node {\small $S_1$};
  \draw[very thick] (1.5,0) ++(90:1) -- ++(90:1) -- ++(60:0.5) ++(60:0.3) node {\small $T_3$};
  \draw[very thick] (1.5,0) ++(0:1) -- ++(0:1) -- ++(30:0.5) ++(30:0.3) node {\small $P_4$};
  \draw[very thick] (1.5,0) ++(0:1) -- ++(0:1) -- ++(330:0.5) ++(330:0.3) node {\small $S_2$};
 
  \draw (0,0.2) node {\small $\chi$};
  \draw (1.5,0) ++ (225:1.3) node {\small $\alpha$};
  \draw (1.5,0) ++ (315:1.3) node {\small $\beta$};
%  \draw (1.5,0) ++ (10:1.15) node {\small $\gamma$};
  \draw (1.5,0) ++ (45:1.3) node {\small $\delta$};
  \draw (1.5,0) ++ (135:1.3) node {\small $\epsilon$};
    \draw (-1.75,0) ++ (315:1.6) node {\small $a$};
    \draw (-1.75,0) ++ (225:1.6) node {\small $b$};
  %  \draw (-1.5,0) ++ (175:1.15) node {\small $c$};
    \draw (-1.75,0) ++ (135:1.6) node {\small $d$};
    \draw (-1.75,0) ++ (45:1.6) node {\small $e$};
  \end{tikzpicture}
 }

By $|a|$ we mean the length of edge $a$, and we use similar notation for the other edges. By isometry of $\mu$ we obtain 
\begin{align}
\label{align:edges_first_cycle}
 |c|+|d|=\frac{3}{4}A-\log|q_{12}|-\frac{1}{2}A=\frac{1}{4}A-\log|q_{12}|=|a|+|\lambda|
\end{align}

and
\begin{align*}
 |\beta|+|\gamma|-|\epsilon|&=-\frac{1}{2}\log|q_{22}|-\frac{1}{4}B-(-\log|q_{22}|+\log|q_{12}|)+(\frac{3}{4}B-\log|q_{12}|) \\
 &= \frac{1}{2}\log|q_{22}|-2\log|q_{12}|+\frac{1}{2}B 
 =-\log|q_{12}|=|\lambda|
 \end{align*}
where $|c|=|\gamma|=|\lambda|=0$ in the connecting edge case.

We start by looking at $\log|f|$. It suffices to calculate the slopes along $a$ and $\alpha$ denoted by $m_a$ and $m_\alpha$ resp. since the other slopes can be determined using the slope formula (Theorem \ref{thm:poincare}), e.g. the slope along edge $b$ is $m_a-1$. Since $\log|f|$ is piecewise linear by walking around the first cycle we get 
\begin{align*}
 m_\lambda |\lambda|+m_a|a|+(m_a-1)|b|+(m_a-2)(|c|+|d|)+(m_a-1)|e|&=0.
\end{align*}
Rearranging and making use of equation (\ref{align:edges_first_cycle}) leads us to 
\begin{align*}
 m_\lambda|\lambda|+m_a(-\log|q_{11}|-|\lambda|)&=-\log|q_{11}|
\end{align*}

Thus $(m_\lambda-m_a)|\lambda|=(1-m_a)\log|q_{11}|$. If $|\lambda|=0$, i.e. if we are in the connecting edge case, then $m_a=1$. If not then $m_\lambda-m_a=m_\alpha$ by the slope formula. Thus
\begin{align}
\label{align:cycle1}
 m_\alpha|\lambda|=(1-m_a)\log|q_{11}|.
\end{align}
Similarly since
\begin{align*}
m_\lambda |\lambda|+m_\alpha|\alpha|+(m_\alpha-1)(|\beta|+|\gamma|)+m_\alpha|\delta|+(m_\alpha+1)|\epsilon|&=0
\end{align*}
and $|\beta|+|\gamma|=|\epsilon|+|\lambda|$ on the second cycle we obtain 
\begin{align*}
m_\lambda |\lambda|+m_\alpha(-\log|q_{22}|-|\lambda|)&=|\lambda| 
\end{align*}
 Thus 
\begin{align*}
-m_\alpha \log|q_{22}|&=(m_\alpha-m_\lambda+1)|\lambda|.
\end{align*}

If $|\lambda|=0$ then $m_\alpha=0$. If not then $m_\alpha-m_\lambda=-m_a$ thus
\begin{align*}
-m_\alpha \log|q_{22}|&=(1-m_a)|\lambda|.
 \end{align*}
Inserting $m_\alpha=-(1-m_a)\dfrac{|\lambda|}{\log|q_{22}|}$ into (\ref{align:cycle1}) gives:
\begin{align*}
-(1-m_a)\dfrac{|\lambda|^2}{\log|q_{22}|}&=(1-m_a)\log|q_{11}| 
\end{align*}
Since $|\lambda|=-\log|q_{12}|$ and $0<\log|q_{12}|^2<\log|q_{11}|\log|q_{22}|$ this implies $m_a~=~1$ and hence $m_\alpha=0$.

The calculation for $\log|g|$ is similiar to that of $\log|f|$ treating the second cycle first.

Thus we obtain the following slope vectors for $(\log|f|,\log|g|)$:
$$
m_a=(1,0),m_b=(0,-1),m_c=(-1,-1),m_d=(-1,0),m_e=(0,1)
$$
on the first cycle and 
$$
m_\alpha=(0,1),m_\beta=(-1,0),m_\gamma=(-1,-1),m_\delta=(0,-1),m_\epsilon=(1,0)
$$
on the second one. If the skeleton has a shared edge then by the slope formula (Theorem \ref{thm:poincare}) the corresponding slope vector is $m_\lambda=(1,1)$. 
If the skeleton has a connecting edge we obtain $m_\chi=(-1,1)$. 
Hence the tropicalization maps the skeleton of the curve isometrically to its image, i.e. the tropicalization is faithful. See figure \ref{fig:trop_shared} for an illustration. 

}
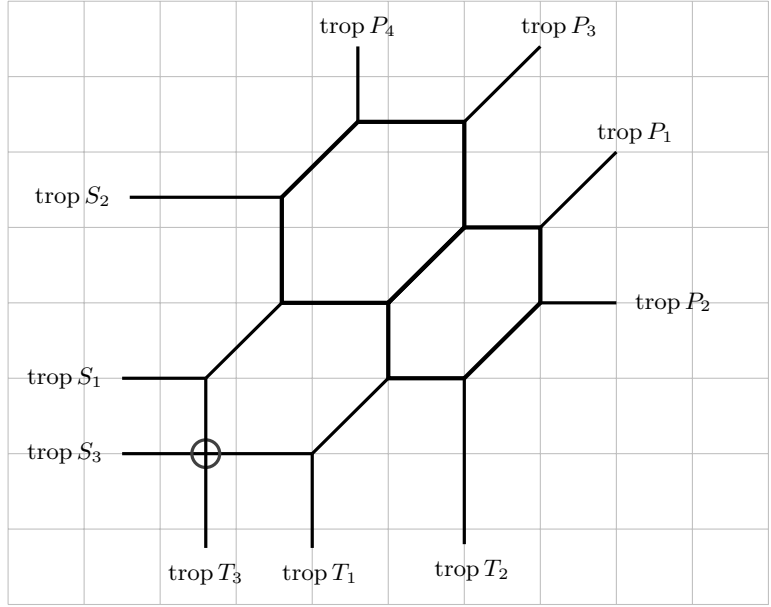
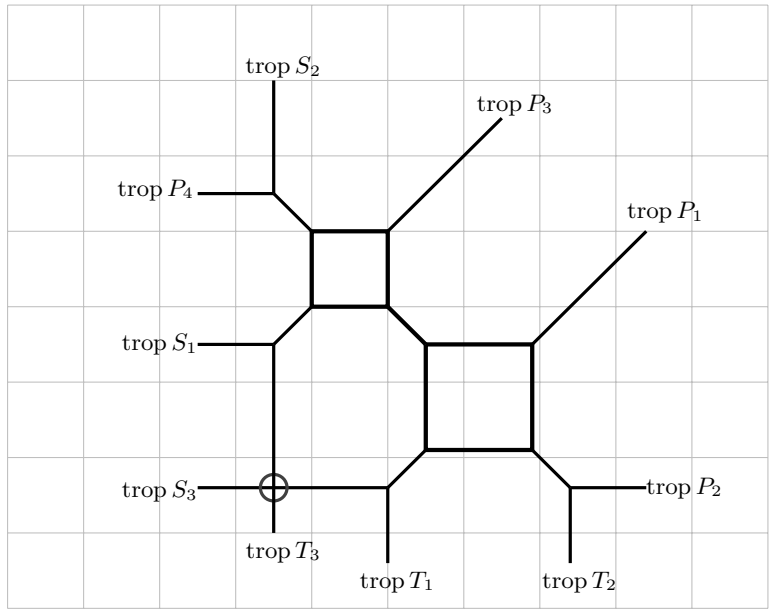
\fig{trop_shared}{Tropicalization with respect to $\log|f|$ and $\log|g|$ of theorem \ref{thm:faithful_skeleton} in the case of a shared edge resp. connecting edge where the possible intersection points are marked by a circle.}{
\begin{tikzpicture}[scale=1]
\draw[lightgray] (-5,-4) grid (5,4);
 \draw[ultra thick] (0,0) -- (1,1) -- (2,1) -- (2,0) -- (1,-1) -- (0,-1) -- (0,0);
 \draw[ultra thick] (0,0) -- (1,1) -- (1,2.4) -- (-0.4,2.4) -- (-1.4,1.4) -- (-1.4,0) -- (0,0);
 \draw[very thick] (2,1) -- (3,2) +(.25,.25) node {\small $\trop P_1$};
 %\draw (2,-1) -- (3,-2);
 \draw[very thick] (1,-1) -- (1,-3.2)+(0.1,-.35) node {\small $\trop T_2$};
 \draw[very thick] (2,0) -- (3,0)+(.75,0) node {\small $\trop P_2$};
 \draw[very thick] (0,-1) -- (-1,-2);
 \draw[very thick] (-1,-2) -- (-1,-3.25)+(0.1,-.35) node {\small $\trop T_1$};
 \draw[very thick] (-1,-2) -- (-3.5,-2)+(-.75,0) node {\small $\trop S_3$};
 \draw[very thick] (1,2.4) -- (2,3.4)+(.25,.25) node {\small $\trop P_3$};
 %\draw (-1,2) -- (-2,3);
 \draw[very thick] (-1.4,1.4) -- (-3.4,1.4) +(-.75,0) node {\small $\trop S_2$};
 \draw[very thick] (-0.4,2.4) -- (-0.4,3.4) +(0,.25) node {\small $\trop P_4$};
 \draw[very thick, darkgray] (-2.4,-2) circle (0.183);
 \draw[very thick] (-1.4,0) -- (-2.4,-1);
 \draw[very thick] (-2.4,-1) -- (-3.5,-1) +(-.75,0) node {\small $\trop S_1$};
 \draw[very thick] (-2.4,-1) -- (-2.4,-3.25) +(0,-.35) node {\small $\trop T_3$};
 \end{tikzpicture}
\\[.5cm]
\begin{tikzpicture}[scale=.5]
 \draw[lightgray,step=2] (-10,-8) grid (10,8);
 \draw[ultra thick] (0,0) -- (0,2) -- (-2,2) -- (-2,0) -- (0,0);
 \draw[ultra thick] (1,-1) -- +(2.8,0) -- +(2.8,-2.8) --+ (0,-2.8) -- (1,-1);
 \draw[ultra thick] (0,0)--(1,-1);
 \draw[very thick] (3.8,-1) -- (6.8,2) +(.5,.5) node {\small $\trop P_1$};
 \draw[very thick] (3.8,-3.8) -- (4.8,-4.8);
 \draw[very thick] (4.8,-4.8) -- (4.8,-6.8)+(.25,-.5) node {\small $\trop T_2$};
 \draw[very thick] (4.8,-4.8) -- (6.8,-4.8)+(1,0) node {\small $\trop P_2$};
 \draw[very thick] (1,-3.8) -- (0,-4.8);
 \draw[very thick] (0,-4.8) -- (0,-6.8)+(.25,-.5) node {\small $\trop T_1$};
 \draw[very thick] (0,-4.8) -- (-5,-4.8)+(-1,-.1) node {\small $\trop S_3$};
 \draw[very thick] (0,2) -- (3,5)+(.35,.35) node {\small $\trop P_3$};
 \draw[very thick] (-2,2) -- (-3,3);
 \draw[very thick] (-3,3) -- (-3,6) +(.25,.35) node {\small $\trop S_2$};
 \draw[very thick] (-3,3) -- (-5,3) +(-1.1,.1) node {\small $\trop P_4$};
 \draw[very thick, darkgray] (-3,-4.8) circle (0.351);
 \draw[very thick] (-2,0) -- (-3,-1);
 \draw[very thick] (-3,-1) -- (-5,-1) +(-1,0) node {\small $\trop S_1$};
 \draw[very thick] (-3,-1) -- (-3,-6) +(.25,-.5) node {\small $\trop T_3$};
 \end{tikzpicture}
}

\Thm{faithful_extended_skeleton}{
Let $D'=D+S_3$. Then there exists a function $h\in L(D')$ such that $(1,f,g,h):\curve~\longhookrightarrow~\PP^3$ is a closed immersion tropicalizing the extended skeleton faithfully. % with respect to $D'$ faithfully.
}
\Prf{
Since $\deg(D')=5$ the divisor is very ample. By Riemann-Roch $\dim L(D')=3$ thus a basis of the linear system defines an embedding into threedimensional projective space. The functions $1,f$ and $g$ are linearly independent in $L(D')$ so they can be completed to a basis by a function having a pole in $S_3$. 
We already know that $\log|f|$ and $\log|g|$ tropicalize the skeleton faithfully. 

First note that there can be up to two additional bounded edges on the extended skeleton with respect to $D'$: The rays to $S_1$ and $T_3$ (resp. $T_1$ and $S_3$) may meet outside the skeleton. The slopes along these edges can be determined using the slope formula thus the slope vectors are $(-1,-1)$ in both cases (oriented away from the skeleton). Hence these edges are mapped isometrically to their image in the tropicalization. 

In the case of a connecting edge there are possibly two more bounded edges between the  point where the rays leading to $S_2$ and $P_4$ meet and the skeleton resp. between the point where the rays leading to $P_2$ and $T_2$ meet and the skeleton. Again using the slope formula we see that the slope vectors are $(-1,1)$ resp. $(1,-1)$ so also these edges are going to be tropicalized faithfully. 

The slope vector of the ray leading to $P_1$ is denoted by $m_{P_1}$ similar notations for the other rays. By convention they are oriented away from the skeleton and are equal to the multiplicities of the poles resp. zeros of each function (Theorem \ref{thm:poincare}). Thus we obtain
$$
m_{P_1}=(1,1),\ m_{P_2}=(1,0),\ m_{T_2}=(0,-1),\ m_{T_1}=(0,-1),\ m_{S_3}=(-1,0)
$$

on the first cycle and
$$
m_{P_3}=(1,1),\ m_{P_4}=(0,1),\ m_{S_2}=(-1,0),\ m_{S_1}=(-1,0),\ m_{T_3}=(0,-1)
$$

on the second cycle.

Without loss of generality choose coordinates in $\RR^2$ such that via the tropicalization map one vertex of the connecting resp. shared edge lies in the origin. Let the plane be divided into two halfspaces by the line $\<(1,1)\>$. Then each cycle lies in one of the two halfspaces. Except for the rays leading to $S_3$ resp. $T_3$ all rays of one cycle lie in the same halfspace. This means these rays cannot intersect with rays from the other cycle. However the rays leading to $S_3$ and $T_3$ may well intersect, see figure  \ref{fig:trop_shared}. By  arguments similar to that in the proof of lemma \ref{lemma:points} there exists a rational function $h$ having poles in $S_3$, $P_1$ and $P_3$ as well as a zero in $T_2$ which implies that none of the other zeros will retract to $\tau(S_3)$. Since $h$ has a pole in $S_3$ the  ray leading to this pole (and also a possible additional bounded edge) will be lifted into the third dimension as the ray leading to $T_3$ will not. %

Thus $\log|h|$ separates the remaining rays, the tropicalization of the extended skeleton with respect to $\log|f|,\log|g|$ and $\log|h|$ is faithful.
}

{\bf Acknowledgements }
This work has been supported by Deutsche Forschungsgemeinschaft (grant WE 4279). I would like to thank Annette Werner for numerous suggestions and corrections.

\end{document}